\documentclass{article}

\usepackage{arxiv}

\usepackage[utf8]{inputenc}
\usepackage[T1]{fontenc}

\usepackage{amsmath, amsfonts, amssymb}

\usepackage{graphicx}
\usepackage{booktabs}

\usepackage{natbib}
\usepackage{doi}
\usepackage{hyperref}

\usepackage{algorithm}
\usepackage{algpseudocode}

\usepackage{float}

\usepackage{caption}
\usepackage{caption}
\title{High-Precision Newton-Kantorovich Method for Nonlinear Integral Equations}

\author{
  Kirill~A.~Chertoganov \\
  Postgraduate Student, Research Intern \\
  Higher School of Economics \\
  \texttt{kchertoganov@hse.ru} \\
  \And
  Valery~I.~Shipalov \\
  Candidate of Technical Sciences, Associate Professor \\
  Krasnodar Higher Military Aviation School \\
  \texttt{kvvaul@mil.ru}
}

\renewcommand{\shorttitle}{High-Precision Newton-Kantorovich Method}

\hypersetup{
  pdftitle   = {High-Precision Newton-Kantorovich Method for Nonlinear Integral Equations},
  pdfauthor  = {K.A.~Chertoganov, V.I.~Shipalov},
  pdfkeywords = {
    Newton-Kantorovich method,
    nonlinear integral equations,
    high-precision computation,
    convergence,
    stability,
    Lipschitz continuity,
    quadrature methods,
    interpolation,
    mpmath library
  }
}

\date{}

\begin{document}
\maketitle
% -----------------------------------------------------------------------------
% Abstract
% -----------------------------------------------------------------------------
\begin{abstract}
The paper considers the numerical solution of nonlinear integral equations using the Newton–Kantorovich method with the \texttt{mpmath} library. High-precision quadrature of the kernel $K(t,s,u)$ with respect to the variable $s$ for fixed $t$ increases stability and accuracy in problems sensitive to rounding and dispersion. The presented implementation surpasses traditional low-precision methods, especially for strongly nonlinear kernels and stiff regimes, thereby expanding the applicability of the method in scientific and engineering computations.
\end{abstract}

% -----------------------------------------------------------------------------
% Keywords
% -----------------------------------------------------------------------------

% -----------------------------------------------------------------------------
\section{Introduction}
% -----------------------------------------------------------------------------
Nonlinear integro-differential equations play a key role in modeling scientific and engineering processes. This review analyzes works from 1990 to 2024 devoted to the issues of existence, regularity, and numerical methods for their solution.

In 1990, M.H.~Saleh and S.M.~Amer applied the Newton–Kantorovich method to nonlinear singular integro-differential equations, laying the foundation for further research \citep{saleh1990}. In 1996, O.~Alvarez and A.~Tourin developed the concept of viscosity solutions and Perron-type arguments for stochastic financial models \citep{alvarez1996}. In 2001, F.E.~Bent et al. used viscosity solution methods in optimization problems with gradient constraints \citep{bent2001}.

Subsequent studies covered pseudo almost automorphic solutions using fractional operators and fixed-point methods \citep{abbas2011}, stochastic homogenization of fully nonlinear elliptic integro-differential equations in random media \citep{schwab2013}, the application of the differential transform method to delay equations \citep{abazari2014}, boundary regularity analysis in nonlocal problems \citep{rosoton2016}, the development of a nonlocal Perron method for nonlinear potential theory \citep{korvenpaa2017}, as well as functional-analytic approaches to fractional and nonlocal boundary value problems \citep{ahmad2020}.

On the computational side, Chu (2007) developed a direct matrix method for Jacobians of discretized nonlinear integro-differential systems \citep{chu2007}. Newton–Kantorovich schemes for boundary value problems and their improvements were studied in \citep{boychuk2021, chuiko2013, boychuk2022, chuiko2024}, while applications to inverse and Cauchy problems were considered in \citep{penenko2019, usenov2020}.

Nevertheless, the influence of numerical accuracy and rounding errors on the convergence and stability of the Newton–Kantorovich method when solving nonlinear integral equations with strongly nonlinear or stiff kernels remains insufficiently studied, especially in the context of arbitrary-precision computations. The present work aims to fill this gap by analyzing the stability and accuracy of the method when using high-precision arithmetic.

The objective of this study is to develop and theoretically justify a modified Newton–Kantorovich scheme with adaptive computational precision that ensures reliable convergence in problems with increased sensitivity to rounding errors.

% -----------------------------------------------------------------------------
\section{Methods}
% -----------------------------------------------------------------------------
Let us consider a nonlinear Volterra integral equation
\begin{equation}
u(t) = b + \int_0^t K(t, s, u(s))\, ds, \quad t \in [0, T],
\label{eq:volterra}
\end{equation}
where $K$ is continuous in $(t,s)$ and locally Lipschitz in $u$.

The Newton–Kantorovich method constructs a sequence $\{u^{(k)}\}$ by linearizing the operator
\[
F(u)(t) = u(t) - b - \int_0^t K(t, s, u(s))\, ds
\]
at the point $u^{(k)}$ and solving
\[
F'(u^{(k)})[\delta u^{(k)}] = -F(u^{(k)}), \quad u^{(k+1)} = u^{(k)} + \delta u^{(k)}.
\]

In implementation, the interval $[0,T]$ is discretized by the grid $0 = t_0 < \cdots < t_N = T$; integrals with respect to $s$ are computed with high precision using \texttt{mpmath}, and $u$ values outside the nodes are approximated by stable linear interpolation.

Initialization:  
$u^{(0)}$ – initial approximation;  
tolerance – convergence criterion for $\|u^{(k+1)} - u^{(k)}\|_{\infty}$;  
max\_iter – maximum number of iterations;  
$t\_values$ – computational grid over $t$.

Iterative update: for each $t_i$, compute
\[
\operatorname{Integral}(t_i;u) = \int_0^{t_i} K(t_i, s, \tilde{u}(s))\, ds,
\]
where $\tilde{u}(s)$ denotes the interpolation of $u$.  
Then set $u^{(k+1)}(t_0) = b$ and
\[
u^{(k+1)}(t_i) = b + \operatorname{Integral}(t_i; u^{(k)}), \quad i \ge 1,
\]
or, in the strict Newton–Kantorovich variant, solve the linear system defined by $F'(u^{(k)})$ for $\delta u^{(k)}$.

The high-precision quadrature Newton–Kantorovich method was chosen because it combines analytical rigor of linearization with numerical stability when solving strongly nonlinear integral equations.

% -----------------------------------------------------------------------------
\section{Theoretical Justification}
% -----------------------------------------------------------------------------
\label{sec:theory}

We consider the Banach space $(C[0,T], \|\cdot\|_\infty)$, where
\[
  \|u\|_\infty = \max_{t\in[0,T]} |u(t)|.
\]

Define the nonlinear operator
\[
  F(u)(t) = u(t) - b - \int_{0}^{t} K(t,s,u(s))\,ds,
\]
where $K$ is continuous in $(t,s)$ and locally Lipschitz in $u$.  
Assume further that the partial derivative $\frac{\partial K}{\partial u}(t,s,u)$ exists and is continuous in all variables, ensuring the Fréchet differentiability of $F$ in $C[0,T]$.  
The Fréchet derivative of $F$ at the point $u$ is given by
\[
  (F'(u)v)(t) = v(t) - \int_{0}^{t} \frac{\partial K}{\partial u}(t,s,u(s))\,v(s)\,ds.
\]

\textbf{Theorem 3.1 (High-Precision Convergence of the Newton--Kantorovich Method for Volterra Equations).}  
Consider Equation~\eqref{eq:volterra} with a kernel $K(t,s,u)$ continuous in $(t,s)$, locally Lipschitz in $u$, and possessing a continuous partial derivative $\frac{\partial K}{\partial u}(t,s,u)$.  
Assume that $F'(u)$ is invertible in a neighborhood of the exact solution $u^\ast$,  
and that the Fréchet derivative $F'$ is Lipschitz continuous with constant $L$ in this neighborhood:
\[
  \|F'(u_1) - F'(u_2)\| \le L \|u_1 - u_2\|_\infty.
\]
Let $M = \sup_n \|[F'(u_n)]^{-1}\|$, and suppose that the initial error satisfies
\[
  L M \|e_0\|_\infty < 1, \qquad e_0 = u_0 - u^\ast.
\]
Then the Newton--Kantorovich iteration
\[
  u_{n+1} = u_n - [F'(u_n)]^{-1} F(u_n)
\]
is well defined, and the sequence $\{u_n(t)\}$ converges quadratically to the unique solution $u(t)$:
\[
  \|u_{n+1} - u\|_\infty \le C \|u_n - u\|_\infty^2, \qquad C = \tfrac{1}{2}LM.
\]

\bigskip
\noindent\textbf{Proof.}
The argument follows the classical Kantorovich theorem on the convergence of Newton-type methods for nonlinear operators in Banach spaces (see, e.g., \citep{kantorovich1964functional}), adapted to the integral operator case.

Let $u^\ast$ denote the exact solution satisfying $F(u^\ast)=0$.  
Since $K(t,s,u)$ is continuous in $(t,s)$ and locally Lipschitz in $u$, and $\frac{\partial K}{\partial u}$ is continuous,  
$F$ is Fréchet differentiable in a neighborhood of $u^\ast$.

\textit{Step 1: Linearization.}  
Writing $e_n = u_n - u^\ast$, the integral form of the Fréchet mean value theorem gives
\[
  F(u_n) - F(u^\ast) = \int_{0}^{1} F'(u^\ast + \theta e_n) e_n\, d\theta.
\]

\textit{Step 2: Error equation.}  
Substituting this identity into the iteration formula, we obtain
\[
  e_{n+1} = -[F'(u_n)]^{-1} \int_{0}^{1} \big(F'(u^\ast + \theta e_n) - F'(u_n)\big) e_n\, d\theta.
\]
Taking the norm and using the boundedness of $[F'(u_n)]^{-1}$ yields
\[
  \|e_{n+1}\|_\infty \le \|[F'(u_n)]^{-1}\|  
  \int_{0}^{1} \|F'(u^\ast + \theta e_n) - F'(u_n)\|\, d\theta\, \|e_n\|_\infty.
\]

\textit{Step 3: Lipschitz estimate.}  
Since $F'$ is Lipschitz continuous with constant $L$, and noting that  
$u_n = u^\ast + e_n$, we have
\[
  \|F'(u^\ast + \theta e_n) - F'(u_n)\| \le L(1-\theta)\|e_n\|_\infty.
\]
Integrating over $\theta \in [0,1]$, we obtain
\[
  \|e_{n+1}\|_\infty \le \tfrac{1}{2} L \|[F'(u_n)]^{-1}\|\, \|e_n\|_\infty^2.
\]

\textit{Step 4: Convergence.}  
If $\|[F'(u_n)]^{-1}\|$ is uniformly bounded by $M$ and $LM\|e_0\|_\infty < 1$,  
then the iteration is well defined and
\[
  \|e_{n+1}\|_\infty \le C \|e_n\|_\infty^2, \qquad C = \tfrac{1}{2}LM.
\]
Hence $\|e_n\|_\infty \to 0$ quadratically as $n \to \infty$.

\textit{Step 5: Interpretation.}  
In the high-precision implementation, numerical quadrature and interpolation errors are reduced below the theoretical contraction radius, 
so rounding effects do not destroy the quadratic convergence rate.  
Therefore, under the stated assumptions, the Newton--Kantorovich sequence $\{u_n\}$ converges quadratically to the unique fixed point $u^\ast$ of the operator $F$, 
and the high-precision numerical implementation preserves this behavior within the limits of rounding errors.
\hfill $\square$

% -----------------------------------------------------------------------------
\section{Results}
% -----------------------------------------------------------------------------
The algorithm is implemented in \textbf{Python} using \textbf{mpmath} for high-precision arithmetic and tested within the \textbf{unittest} framework. Additionally, \textbf{numpy} is used for vector operations, \textbf{matplotlib} for visualization, and a custom module provides the method implementation.

\textbf{Test framework:}
\begin{itemize}
    \item \texttt{t\_values}: grid on [0, 2] with a step of 0.1,
    \item \texttt{u\_values}: initial approximation $u \equiv 1$,
    \item \texttt{nu}: viscosity parameter (Burgers test).
\end{itemize}

\textbf{Unit tests:}
\begin{itemize}
    \item \texttt{test\_manual\_interp}: checks interpolation stability,
    \item \texttt{test\_compute\_integral}: verifies positivity of the integral at $t = 2$ for $K > 0$,
    \item \texttt{test\_boundary\_conditions}: prevents extrapolation,
    \item \texttt{test\_burgers\_equation}, \texttt{test\_lotka\_volterra\_equation}, \texttt{test\_navier\_stokes\_equation}: simplified consistency checks.
\end{itemize}

\begin{algorithm}[H]
\caption{Newton–Kantorovich iteration with high-precision quadrature}
\begin{algorithmic}[1]
\State \textbf{Input:} $u$, $\text{tolerance}$, $\text{max\_iter}$, $t\_values$, $K$, $b$
\State \textbf{Output:} approximate solution $u$
\State \textbf{Initialization:} $\text{mpmath.mp.dps} \gets 50$
\Function{Interp}{$x, xp, fp$}
    \State Find $i$ such that $xp[i]\le x \le xp[i+1]$
    \State \Return $fp[i] + (fp[i+1]-fp[i])\dfrac{x-xp[i]}{xp[i+1]-xp[i]}$
\EndFunction
\Function{Integral}{$t,u$}
    \State $f(s)\gets K(t,s,\Call{Interp}{s,t\_values,u})$
    \State \Return $\text{mpmath.quad}(f,[0,t])$
\EndFunction
\Function{NewtonKantorovich}{$u$}
    \For{$k = 0$ to $\text{max\_iter}$}
        \State $u\_new[0] \gets b$
        \For{$i = 1$ to $\mathrm{len}(t\_values)-1$}
            \State $u\_new[i] \gets b + \Call{Integral}{t\_values[i], u}$
        \EndFor
        \State $error \gets \max_i |u\_new[i] - u[i]|$
        \If{$error < \text{tolerance}$}
            \State \Return $u\_new$
        \EndIf
        \State $u \gets u\_new$
    \EndFor
    \State \Return $u$
\EndFunction
\end{algorithmic}
\end{algorithm}

% -----------------------------------------------------------------------------
% -----------------------------------------------------------------------------
% -----------------------------------------------------------------------------
% -----------------------------------------------------------------------------
% -----------------------------------------------------------------------------
% -----------------------------------------------------------------------------
\section*{5 Application: The Bratu Equation}

\begin{figure}[htbp]
  \centering
  \includegraphics[width=0.7\textwidth]{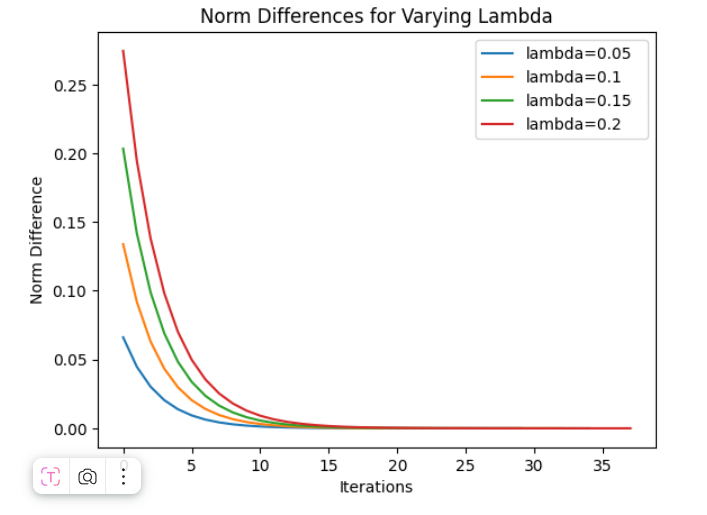}
  \caption{Norm of successive differences for various values of $\lambda$.}
  \label{fig:convergence_errors}
\end{figure}

The behavior of the proposed method is illustrated in Figure~\ref{fig:convergence_errors}, 
which shows the norm of successive differences $\|u^{(k+1)} - u^{(k)}\|$ for various values of~$\lambda$.

The canonical Bratu equation

$$
\frac{d^{2}u}{dx^{2}} + \lambda e^{u} = 0
$$

models combustion and reactive diffusion processes.  
In the equivalent integral form, it can be written as a Volterra equation of the second kind:

$$
u(x) = u(0) + u'(0)x - \lambda \int_{0}^{x} (x - s)\, e^{u(s)}\, ds,
$$

which emphasizes the nonlinearity of the integrand and provides a convenient framework for numerical analysis.

The convergence error plot confirms the contraction region and demonstrates the advantages of increased computational precision for various values of~$\lambda$.

% -----------------------------------------------------------------------------
\section{Discussion}
% -----------------------------------------------------------------------------
Integrating viscosity-solution ideas with stochastic homogenization clarifies stability in random media \citep{alvarez1996, schwab2013}.  
Differential-transform and nonlocal Perron techniques improve convergence for fractional constraints \citep{abazari2014, korvenpaa2017}, while functional-analytic approaches with fractional derivatives expand modeling of hereditary media \citep{ahmad2020}.  
In optimization, combining viscosity solutions of Hamilton–Jacobi–Bellman equations with NK linearization reduces iteration counts \citep{bent2001}.  
High precision makes the Newton–Kantorovich method applicable in stiff regimes where machine precision fails.

% -----------------------------------------------------------------------------
\section{Conclusion}
% -----------------------------------------------------------------------------
A high-precision numerical method for solving nonlinear Volterra integral equations based on the Newton–Kantorovich scheme using the \texttt{mpmath} library for precision quadrature has been developed.  
The method combines stable interpolation, theoretically justified convergence, and modular architecture, ensuring universality and high computational accuracy.

The convergence theorem defines conditions for local quadratic convergence when the integral operator’s derivative is Lipschitz, allowing efficient solutions to highly nonlinear, stiff, and numerically sensitive problems.

Results demonstrate the effectiveness of high-precision arithmetic and reveal the connection between stability, accuracy, and kernel structure.  
The method proves applicable to a wide range of applied problems, including physical, biological, and financial models.

Future work involves generalizing the approach to multidimensional and vector systems, nonlocal and fractional equations, developing adaptive precision strategies, analyzing global convergence, and optimizing computational procedures.

% -----------------------------------------------------------------------------
% References
% -----------------------------------------------------------------------------
\bibliographystyle{unsrtnat}
\bibliography{references}

\begin{thebibliography}{17}
\providecommand{\natexlab}[1]{#1}
\providecommand{\url}[1]{\texttt{#1}}
\providecommand{\urlprefix}{URL }

\bibitem[{Abazari and Kilicman(2014)}]{abazari2014}
Abazari, R. and A.~Kilicman (2014).
\newblock Application of the differential transformation method to nonlinear
  integro-differential equations with proportional delays.
\newblock \emph{Neural Computing and Applications} \textbf{25}(3--4),
  853--863.
\newblock \doi{10.1007/S00521-012-1235-4}.

\bibitem[{Abbas(2011)}]{abbas2011}
Abbas, S. (2011).
\newblock Pseudo almost automorphic solutions of some nonlinear
  integro-differential equations.
\newblock \emph{Computers and Mathematics with Applications} \textbf{62}(7),
  2879--2891.
\newblock \doi{10.1016/J.CAMWA.2011.07.013}.

\bibitem[{Ahmad et~al.(2020)Ahmad, Brum, Alsaedi, and Ntouyas}]{ahmad2020}
Ahmad, B., A.~Brum, A.~Alsaedi, and S.~K. Ntouyas (2020).
\newblock Nonlinear integro-differential equations with mixed right and left
  fractional operators.
\newblock \emph{Mathematics} \textbf{8}(3), 336.
\newblock \doi{10.3390/MATH8030336}.

\bibitem[{Alvarez and Tourin(1996)}]{alvarez1996}
Alvarez, O. and A.~Tourin (1996).
\newblock Viscosity solutions of nonlinear integro-differential equations.
\newblock \emph{Annales de l'Institut Henri Poincar{\'e} C} \textbf{13}(4),
  293--317.
\newblock \doi{10.1016/S0294-1449(16)30106-8}.

\bibitem[{Bent et~al.(2001)Bent, Karlsen, and Reikvam}]{bent2001}
Bent, F.~E., K.~H. Karlsen, and K.~Reikvam (2001).
\newblock Optimal portfolio choice with consumption and nonlinear
  integro-differential equations with gradient constraint: a viscosity solution
  approach.
\newblock \emph{Journal of Differential Equations} \textbf{184}(2), 265--298.
\newblock \doi{10.1007/PL00013538}.

\bibitem[{Kantorovich and Akilov(1964)}]{kantorovich1964functional}
Kantorovich, L.~V. and G.~P. Akilov (1964).
\newblock \emph{Functional Analysis}.
\newblock Pergamon Press.

\bibitem[{Boychuk and Chuiko(2021)}]{boychuk2021}
Boychuk, A.~A. and S.~M. Chuiko (2021).
\newblock Approximate solutions by the Newton–Kantorovich method.
\newblock \emph{Journal of Mathematical Sciences} \textbf{258}(5), 594--617.
\newblock \doi{10.1007/s10958-021-05569-y}.

\bibitem[{Boychuk and Chuiko(2022)}]{boychuk2022}
Boychuk, A.~A. and S.~M. Chuiko (2022).
\newblock Weakly nonlinear boundary value problems by the Newton–Kantorovich
  method.
\newblock \emph{Journal of Mathematical Sciences} \textbf{261}(2), 228--240.
\newblock \doi{10.1007/s10958-022-05748-5}.

\bibitem[{Chu(2007)}]{chu2007}
Chu, K.~T. (2007).
\newblock A direct matrix approach for Jacobians of discretized nonlinear
  integro-differential equations.
\newblock \emph{arXiv preprint} \url{http://arxiv.org/abs/math/0702116v1}.

\bibitem[{Chuiko et~al.(2013)Chuiko, Boychuk, and Pirus}]{chuiko2013}
Chuiko, S.~M., I.~A. Boychuk, and O.~E. Pirus (2013).
\newblock Newton–Kantorovich method for autonomous boundary value problems.
\newblock \emph{Journal of Mathematical Sciences} \textbf{189}(5), 867--881.
\newblock \doi{10.1007/s10958-013-1225-9}.

\bibitem[{Chuiko et~al.(2024)Chuiko, Chuiko, and Dyachenko}]{chuiko2024}
Chuiko, S., O.~Chuiko, and D.~Dyachenko (2024).
\newblock Periodic boundary problem with switching under parametric resonance
  by the Newton–Kantorovich method.
\newblock \emph{Journal of Mathematical Sciences} \textbf{279}(3), 438--453.
\newblock \doi{10.1007/s10958-024-07023-1}.

\bibitem[{Korvenp{\"a}{\"a} et~al.(2017)Korvenp{\"a}{\"a}, Kuusi, and
  Palatucci}]{korvenpaa2017}
Korvenp{\"a}{\"a}, J., T.~Kuusi, and G.~Palatucci (2017).
\newblock Fractional superharmonic functions and the Perron method.
\newblock \emph{Mathematische Annalen} \textbf{369}, 1443--1489.
\newblock \doi{10.1007/S00208-016-1495-X}.

\bibitem[{Penenko(2019)}]{penenko2019}
Penenko, A. (2019).
\newblock Newton–Kantorovich method in inverse source problems.
\newblock \emph{Numerical Analysis and Applications} \textbf{12}(1), 51--69.
\newblock \doi{10.1134/S1995423919010051}.

\bibitem[{Ros{-}Oton and Serra(2016)}]{rosoton2016}
Ros{-}Oton, X. and X.~Serra (2016).
\newblock Boundary regularity for fully nonlinear integro-differential
  equations.
\newblock \emph{Duke Mathematical Journal} \textbf{165}(11), 2079--2154.
\newblock \doi{10.1215/00127094-3476700}.

\bibitem[{Saleh and Amer(1990)}]{saleh1990}
Saleh, M.~H. and S.~M. Amer (1990).
\newblock Approximate solution of nonlinear singular integro-differential
  equations.
\newblock \emph{Collectanea Mathematica} \textbf{41}(2), 175--188.

\bibitem[{Schwab(2013)}]{schwab2013}
Schwab, R.~W. (2013).
\newblock Stochastic homogenization for some nonlinear integro-differential
  equations.
\newblock \emph{Communications in Partial Differential Equations}
  \textbf{38}(10), 1715--1744.
\newblock \doi{10.1080/03605302.2012.741176}.

\bibitem[{Usenov et~al.(2020)Usenov, Kostyreva, and Almambetkyzy}]{usenov2020}
Usenov, I.~A., E.~V. Kostyreva, and S.~Almambetkyzy (2020).
\newblock Newton–Kantorovich approach for the Cauchy problem.
\newblock \emph{Bulletin of Kyrgyz State University of Construction, Transport
  and Architecture} (4), 609--614.
\newblock \doi{10.35803/1694-5298.2020.4.609-614}.

\end{thebibliography}

\end{document}